\documentclass{amsart}

\usepackage{amsmath}
\usepackage{amsthm}

\newcommand{\Atlas}{{$\mathbb{ATLAS}$}}
\newcommand{\GAP}{\textsf{GAP}}
\newcommand{\MAGMA}{\textsf{MAGMA}}
\newcommand{\B}{{\mathbb B}}
\newcommand{\M}{{\mathbb M}}
\newcommand{\C}{{\mathbb C}}

\def \ov {\overline}

\def \udot {{}^{\textstyle .}}

\newcommand{\Co}{\mathrm{Co}}

\newcommand{\HN}{\mathrm {HN}}
\newcommand{\HS}{\mathrm {HS}}
\newcommand{\Fi}{\mathrm {Fi}}

\newcommand{\GL}{\mathrm{GL}}
\newcommand{\Sp}{\mathrm{S}}
\newcommand{\Uni}{\mathrm{U}}
\newcommand{\Lin}{\mathrm{L}}

\title{Verification of the Ordinary Character Table of the Baby Monster}
\author{Thomas Breuer \and Kay Magaard 
\and Robert A. Wilson}

\begin{document}
\begin{abstract}
We prove the correctness of the character table of the sporadic simple
Baby Monster group that is shown in the {\Atlas} of Finite Groups.

\end{abstract}

\maketitle

\begin{quote}
\textit{
In memory of
our friend and colleague Kay Magaard, who sadly passed away 
during the preparation of this paper.}
\end{quote}

\section{Introduction}

Jean-Pierre Serre has raised the question of verification
of the ordinary character tables
that are shown in the {\Atlas} of Finite Groups \cite{Atlas}.
This question was partially answered in the paper \cite{BMO17},
the remaining open cases being the largest two sporadic simple groups,
the Baby Monster group $\B$ and the Monster Group $\M$,
and the double cover $2.\B$ of $\B$.

The current paper 
describes a verification of the character table of $\B$.
The computations shown in \cite{ctblatlas} then imply that also the
{\Atlas} character table of $2.\B$ is correct.
As in \cite{BMO17},
one of our aims is to provide the necessary data in a way that makes it
easy to reproduce our computations.

The {\Atlas} character table of the Baby Monster derives from the 
original calculation of the conjugacy classes and rational character
table by David Hunt, described very briefly in \cite{hunt}. The irrationalities
were calculated by 
the CAS team in Aachen \cite{Neubueser}.

\section{Strategy}

We begin with a preliminary section, Section \ref{sect:verifypresentation},
whose aim is to prove that certain specified matrices do indeed
generate copies of the Baby Monster. These matrices can then be used
in the main computation.
The $Y_{555}$ presentation of the BiMonster
implies a $Y_{433}$ presentation for $\B$ (see~\cite{Ivanov}),
given that the Schur multiplier $H^2(\M,\C^*)$ of the Monster has odd order.
The Schur multiplier of the Monster was calculated by Griess \cite{schurmult}.
We use the $Y_{433}$ presentation to prove that three pairs of matrices,
of dimension $4370$ over the field with two elements
and of dimension $4371$ over the fields with three and five elements,
respectively, generate the group $\B$, and, moreover, that mapping one pair of these
generators to any other such pair defines a group isomorphism.

In Section~\ref{sect:classes} 
we compute a first approximation to the list of conjugacy class names, by
establishing invariants, in terms of the
above three matrix representations of $\B$,
that in fact distinguish almost all conjugacy classes of cyclic subgroups of $\B$.
These invariants are then used to 
determine the power maps between the specified unions of conjugacy classes.

In order to compute the conjugacy classes of $\B$ and the
corresponding centralizer orders,
we apply the following general statement.
For a group $G$ and an element $g \in G$,
if $x, y \in G$ power to $g$ then $x \sim_G y$ if and only if
$x \sim_N y$ where $N:= N_G(\langle g \rangle)$;
moreover $C_G(x) = C_N(x)$.
This implies that it suffices to find the normalizers (or overgroups thereof)
of prime order subgroups, and their character tables.
Section~\ref{sect:oddcentralizers} deals with 
the first problem, 
Section~\ref{classlist}
with the second.

At this point,
we know the conjugacy classes of $\B$, and their lengths.
Some further calculations then match these up with the
names listed in Section~\ref{sect:classes}, by using the
given invariants and some small extra arguments.
Finally, the irreducible characters of $\B$ 
are computed in Section~\ref{sect:irreducibles},
using character theoretic methods such as induction from several subgroups
of $\B$.

\section{Verifying a presentation for the Baby Monster}%
\label{sect:verifypresentation}
In this section
we give words in the `standard generators' for the Baby Monster,
that represent the $11$ transpositions in the $Y_{433}$ presentation.
This provides a relatively straightforward test to prove that 
a given black-box group is in fact isomorphic to the Baby Monster.

\subsection{The presentation}
A presentation for the Baby Monster sporadic simple group $\B$
was conjectured in the {\Atlas} \cite{Atlas}, 
and proved by Ivanov \cite{Ivanov}, subject to the Monster not having a proper
double cover. This hypothesis has been proved by Griess \cite{schurmult}. 
 
The presentation is on $11$ generators $t_i$ ($1\le i\le 11$), satisfying
the Coxeter relations $t_i^2=1$ for all $i$, $(t_it_j)^3=1$ for $(i,j)=(1,2)$,
$(2,3)$, $(3,4)$, $(4,5)$, $(5,6)$, $(6,7)$, $(7,8)$, $(5,9)$, $(9,10)$, $(10,11)$,
and $(t_it_j)^2=1$ for $i<j$ otherwise. Adjoining one extra relation, 
$(t_5t_4t_3t_5t_6t_7t_5t_9t_{10})^{10}=1$,
nicknamed the
`spider relation', gives a presentation for $2\times 2\udot \B$.
To obtain a presentation for $\B$ itself,
we need two extra relations, $(t_5t_4t_3t_6t_7t_8t_9)^9=1$
and $(t_5t_4t_3t_6t_9t_{10}t_{11})^9=1$.
Since the Coxeter diagram has three `arms', of lengths $4,3,3$, this
presentation is known as the $Y_{433}$ presentation.

Matrices generating a copy of the Baby Monster were first produced in the
early 1990s \cite{Bmod2}. These act on a vector space of dimension $4370$
over the field of order $2$. In order to prove, without relying on
the character table, that these matrices do indeed generate the Baby Monster,
a method was given for producing elements
of this group that satisfy the $Y_{433}$ presentation for the Baby Monster.
However, actual words for these elements
were not given in \cite{Bmod2}.

In this section 
we rectify this deficiency in \cite{Bmod2}, and hence enable the reader to check
relatively easily that the matrices given in \cite{webatlas}, that are
claimed to generate the Baby Monster in various different representations,
do in fact generate the Baby Monster. In addition to the representation over
the field of order $2$, already mentioned, we checked the representations
over the fields of order $3$ and $5$ constructed in \cite{Bmod35}.
All three of these representations will be required later on, for determining
certain class fusions and power maps.

We begin with the `standard generators' in the sense of \cite{stg,webatlas},
that is an element $a\in 2C$ and an element $b\in 3A$ such that $ab$ has 
order $55$ and $(ab)^4bab(ab^2)^2$ has order $23$. The cited references
explain how to find such generators in a group which is in fact
isomorphic to the Baby Monster. All calculations described in
this section were performed using the
C Meataxe written by Michael Ringe \cite{cmtx}, based on the original
Meataxe of Richard Parker \cite{mtx}.

\subsection{Finding the generators for the presentation}
The calculations in this section were performed using the standard generators for
(a group that is claimed to be)
the Baby Monster in its $4370$-dimensional representation over the field of order $2$,
taken from \cite{webatlas}. 
Following the 10-step method described in \cite{Bmod2}, we proceed as follows.
Steps 1--4 are devoted to finding generators $t_1,\ldots,t_8,t_{11}$ for
a particular subgroup $2\times S_9$.
In Steps 5--8 we centralize successively the elements $t_1$, $t_3t_4$, $t_6t_7$
and $t_8$ to produce a small number of candidates for $t_9$ and $t_{10}$.
These candidates are tested in Step 9, at which point all the required generators
have been found. 
Step 10 tests the relations.

\subsubsection*{Step 1}

\textit{Take an arbitrary $2A$-element, and call it $t_{11}$. Find
$C(t_{11})\cong 2\udot{}^2E_6(2){:}2$.}

The element $d=(ab)^{15}b$ has order $38$, and powers to the involution
$t_{11}=d^{19}$. The centralizer of $t_{11}$ is generated by $d$ and $c=(at_{11})^3$.

\subsubsection*{Step 2}

\textit{Find a subgroup $F\cong Fi_{22}{:}2$ inside $C(t_{11})$.}

We restrict the representation of $\B$ to $H:=\langle c,d\rangle\cong 2\udot {}^2E_6(2){:}2$,
find the composition factors using the \texttt{chop} program of the Meataxe,
and extract a $78$-dimensional irreducible representation of the quotient
$\ov{H}:={}^2E_6(2){:}2$, in which the computational searches for steps 2--4 are performed.
(This use of a small representation reduces the computation time by a factor of around $10^5$.) 
The invariant $q(x)=\mathrm{rank}(1+x)$ is useful for identifying
conjugacy classes. In particular, $\ov{cd^5}$ is an element of order $26$ in the outer half of $\ov{H}$,
so powers to an involution $x\in 2D$, in the {\Atlas} notation
for conjugacy classes in ${}^2E_6(2){:}2$. We calculate $q(x)=26$. As $\ov{c}$ is an involution
in the outer half and $q(\ov{c})=36$, we deduce that $\ov{c}\in 2E$.

Similarly, $\ov{cd}$ has order $36$ and therefore powers to a $3C$-element $x$, with $q(x)=54$,
while $\ov{cd^2}$ has order $48$ and therefore powers to a $3A$-element $y$, with $q(y)=42$.
We then find that $\ov{cd^3}$ has order $30$ and powers to an element $z$ of order $3$
with $q(z)=48$. Hence $z$ is in class $3B$. 
Looking at a few groups generated by conjugates of $\ov{c}$ and $\ov{(cd^3)^{10}}$
we quickly find that if $\ov{e}=\ov{((cd^3)^{10})^d}$
then $F:=\langle \ov{c},\ov{e}\rangle\cong \Fi_{22}{:}2$.
\subsubsection*{Step 3}

\textit{Find a subgroup $S\cong S_{10}$ inside $F$.}

The elements $e$ and $cec$ then generate the subgroup $\Fi_{22}$ of index $2$ in $F$,
in which we find  
$\ov{((ec)^6c(ec)^3)^2(ec)^2cec}$
is an element of order $10$ powering to an element of $\Fi_{22}$-class $2A$.
The element $\ov{(ec)^8c(ec)^3}$ has order $9$
and is most likely to be in class $9C$ in $\Fi_{22}$. 
Since there is no simple test for this, we proceed
and hope for the best.
Looking at conjugates of these elements we soon find a pair $f,g$ generating 
$S:=\langle f,g\rangle\cong S_{10}$,
as follows.
\begin{eqnarray*}
\ov{f}&=&  
\ov{(((ec)^6c(ec)^3)^2(ec)^2cec)^{5})^{(ec)^4}}\cr
\ov{g}&=&\ov{((ec)^8c(ec)^3)^{(ece^2c)^2}}
\end{eqnarray*}
\subsubsection*{Step 4} 

\textit{In $\langle \ov{t_{11}}\rangle \times S$ find transpositions $\ov{t_1},\ldots, \ov{t_8}$ generating
$T\cong S_9$ with the required Coxeter relations.}

These transpositions can be taken as
$\ov{t_1}=\ov{f}$ and $\ov{t_{n+2}}=\ov{f^{gfg^n}}$ for $0\le n\le 6$.

\subsubsection*{Step 5}

\textit{Find $C=C(t_1)\cong 2\udot{}^2E_6(2){:}2$.}

 This step has to be carried out in the $4370$ dimensional
representation.
Standard dihedral group methods give the element $p=((ab)^5t_{11}(ab)^{-5}t_1(ab)^5)^{-1}$
that conjugates $t_{11}$ to $t_1$, so that $h=c^p$ and $i=d^p$ generate 
the group $C:=\langle h,i\rangle=C(t_1)$.

\subsubsection*{Step 6}

\textit{Find $D=C_{C'}(t_3t_4)\cong 6 \times U_6(2)$.} 

Again we restrict the representation of $\B$ to $\langle h,i\rangle$, and extract a copy
of the $78$-dimensional representation, in which
to carry out steps 6--8. Following the instructions in \cite{Bmod2} we found
two elements of the centralizer of $\ov{t_3t_4}$ to be
$\ov{j}=\ov{[(t_5)^{i^2},t_3t_4]}$ and $\ov{k}=\ov{[(t_5)^{i^5},t_3t_4]}$.
Together with $\ov{t_6t_7}$ and $\ov{t_8}$, these are enough to generate $\ov{D}\cong3\times U_6(2)$.

\subsubsection*{Step 7}

\textit{Find $E=C_{D'}(t_6t_7)\cong 3\times U_4(2)$.}

 Similarly we found the following elements centralizing $\ov{t_6t_7}$:
\begin{eqnarray*}
\ov{l}&=& \ov{[t_8^{jk},t_6t_7][t_8^{kj},t_6t_7]}\cr
\ov{l_3}&=&\ov{[t_8^{(jk)^4},t_6t_7]}\cr
\ov{l_4}&=&\ov{t_8^{(jk)^3kj}}
\end{eqnarray*}
These are sufficient to generate $\ov{E}\cong 3\times U_4(2)$. 

\subsubsection*{Step 8} 

\textit{Find the twelve [sic] transpositions in $E$ which commute with $t_8$.}

There is a slight error in \cite{Bmod2} at this point. There are
in fact $13$ transpositions in $U_4(2)$ that centralize $t_8$, rather than
$12$ as stated there. They are the
$13$ transpositions in a copy of $2\udot(A_4\times A_4).2$. One is the central
involution and the other $12$ are in the outer half. Together they generate
a subgroup $2^{1+4}{:}S_3$ of index $3$. Presumably the central involution
was omitted from the original calculation. However, it commutes with neither $t_5$ nor $t_{11}$,
so is not a candidate for $t_9$ or $t_{10}$.

First we looked for conjugates of $\ov{l_4}$ that commute with $\ov{t_8}$, using the
elements $\ov{l_5}=\ov{(ll_3l_4)^3l_3l_4}$ and $\ov{l_6}=\ov{(ll_3l_4)^2l_3l_4ll_3l_4}$,
both of order $9$,
for the conjugation. Then the four conjugates $\ov{m_1}=\ov{l_4}$, $\ov{m_2}=\ov{l_4^{l_5^4}}$,
$\ov{m_3}=\ov{l_4^{l_5^5}}$ and $\ov{m_4}=\ov{l_4^{l_5^2l_6^2}}$ are sufficient to generate
$2^{1+4}{:}S_3$.
The central involution is $\ov{(m_1m_2m_3m_4)^3}$. Then $\ov{m_1}$ and its conjugates by
$\ov{m_2m_3}$, $\ov{m_2m_4}$, and $\ov{m_3m_4}$ are four transpositions mapping to the same
involution in the quotient $S_3$, so give the full set of $12$ after conjugating by
$\ov{m_1m_2}$ and $\ov{m_2m_1}$.
 
\subsubsection*{Step 9} 

\textit{Check these twelve [sic] transpositions for candidates for $t_9$ and $t_{10}$.
There is only one possibility up to an obvious inner automorphism.}

This step was again carried out in the $4370$ dimensional representation.
The calculations that do not involve $t_2$ could have been done in $78$ dimensions, but
the time saved would be of the order of one minute, so insignificant.
Of the $13$ transpositions, the only one which commutes
with $t_2$ and $t_{11}$ but not $t_5$ is $m_1^{m_2m_3m_2m_1}$. Hence this is
the only possibility for $t_9$. There are two that commute with $t_5$ but not $t_{11}$,
namely $m_1^{m_2m_3}$ and $m_1^{m_2m_4m_1m_2}$. But $t_{11}$ conjugates one to the
other, so without loss of generality we may take $t_{10}=m_1^{m_2m_3}$.

\subsection{Verifying the presentation}
We now have a straight line program for producing the elements $t_1,\ldots,t_{11}$
from the elements $a,b$. This program is given in Table~\ref{words1}
for convenience. It must be applied in every claimed
representation of the Baby Monster, and then the relations of the presentation
must be checked.

\begin{table}
\caption{\label{words1}Words to express $t_1,\ldots, t_{11}$ in terms of $a,b$}
$$
\begin{array}{rcl}
\hline
d&=&(ab)^{15}b\cr
t_{11}&=&d^{19}\cr
c&=&(at_{11})^3\cr
e&=&((cd^3)^{10})^d\cr
t_1=f&=&(((ec)^6c(ec)^3)^2ece^2c)^{5})^{(ec)^4}\cr
g&=&((ec)^8c(ec)^3)^{(ece^2c)^2}\cr
t_{n+2}&=&f^{gfg^n}\mbox{ for }0\le n\le 6\cr
p&=&{((ab)^5t_{11}(ab)^{-5}t_1(ab)^5)^{-1}}\cr
h&=&c^p\cr
i&=&d^p\cr
j&=&[(t_5)^{i^2},t_3t_4]\cr
k&=&[(t_5)^{i^5},t_3t_4]\cr 
l&=& [(t_8)^{jk},t_6t_7][(t_8)^{kj},t_6t_7]\cr
l_3&=&[(t_8)^{(jk)^4},t_6t_7]\cr
l_4&=&(t_8)^{(jk)^3kj}\cr
l_5&=&(ll_3l_4)^3l_3l_4\cr
m_2&=&(l_4)^{l_5^4}\cr
m_3&=&(m_2)^{l_5}\cr 
t_{10}&=&m_3m_2l_4m_2m_3\cr
t_9&=&l_4m_2t_{10}m_2l_4\cr
\hline
\end{array}
$$\end{table}
\subsubsection*{Step 10} 

\textit{Prove that $t_1,\ldots,t_{11}$ satisfy all the required relations.}

We check the $66$ Coxeter relations by finding the order of the
elements $t_it_j$ for all $j\ge i$. (The relations $t_i^2=1$ are implicit in the calculation,
but were explicitly checked again.) Similarly the spider relation is
checked by confirming that the element $t_5t_4t_3t_5t_6t_7t_5t_9t_{10}$ has order $10$.
Finally, we check that $t_5t_4t_3t_6t_7t_8t_9$ and $t_5t_4t_3t_6t_9t_{10}t_{11}$ have
order $9$. (In the particular representations we checked, this last check can be
omitted, since it is straightforward to show in each case that the centre of
the group is trivial. Indeed, Schur's Lemma implies the centre consists of scalars,
while the generators have determinant $1$ and the only scalar of determinant $1$ is $1$.)

We verified the relations in the three representations from \cite{webatlas},
that is in dimension $4370$ over the field of order $2$, and in dimension
$4371$ over the fields of orders $3$ and $5$. The total computation time was under
$12$ hours.

\subsection{Reversing the process}
To complete the proof that the matrices given in \cite{webatlas}
generate the Baby Monster, we have to reverse the process, and find the
Atlas standard generators in terms of the $Y_{433}$ generators.
First we make two elements $r=t_1t_2t_3t_4t_5t_6t_7t_8$ of
order $9$ and $s=t_5t_9t_{10}t_{11}$ of order $5$.
These elements in fact generate the whole group, but
it is not necessary to prove this at this stage.
Now $t_1t_2$ is an element in class $3A$. A random search produces
the element $(r^7s)^{15}$ in class $2C$. (Again, it is not necessary
to prove that this involution is in class $2C$. However, we used the
conjugacy class invariants in \cite{Bclass} to guide us, and found
that this element $x$ has $q(x)=2158$, which identifies the conjugacy
class as $2C$.)
Another random search gives a candidate pair of standard generators
\begin{eqnarray*}
a'&=& (r^7s)^{15}\cr
b'&=&(t_1t_2)^{(sr)^{10}}
\end{eqnarray*}

Finally, we use the \texttt{chop} program of the Meataxe to conjugate the
matrices to a standard basis with respect to, first, the generators $a,b$,
and then, the generators $a',b'$. This calculation must, of course,
be carried out in each representation that we wish to check.
We found that in all three of the representations in \cite{webatlas},
the resulting pairs of matrices
are identical, proving that all claimed generating sets do indeed
generate the same group. (This does not mean that $a',b'$ are the
same elements as $a,b$, merely that the pair $(a',b')$ is equivalent to
$(a,b)$ under an automorphism of the group, and therefore
under conjugation.) The total computation time was under $3$ hours.

\section{Conjugacy class invariants and power maps 
in the Baby Monster}%
\label{sect:classes}
\subsection{Introduction}
In this section,
we produce a list of easily computed
conjugacy class invariants for
a specified list of 
elements of the Baby Monster,
which are in fact good enough to distinguish all
conjugacy classes of cyclic subgroups except $16D$ and $16F$.
As a result, we have a splitting of the elements into small unions
of conjugacy classes, and power maps between these unions of classes.
The final splitting into conjugacy classes, and refinements of
power maps, is done later.

In \cite{Bclass} (and also in \cite{webatlas})
there is a list of words in the
Atlas standard generators of the Baby Monster, suitable powers of which
are in fact representatives for the $184$ conjugacy classes.
However, the proof given there depends on the accuracy
of the {\Atlas} character table of the Baby Monster, and in particular
on the accuracy of the power map information.
It is therefore necessary to provide a new proof, which does 
not depend on the character table.
We can of course use the \emph{words}, as long as we do not
quote from \cite{Bclass} any of the properties of the corresponding
elements of the Baby Monster.
We 
assume that the three representations
of the Baby Monster given in \cite{webatlas} do indeed
represent the Baby Monster. This was proved in 
Section \ref{sect:verifypresentation}.

\subsection{The words and their names}
In \cite{Bclass} there is a list of $76$ words for elements of specified orders,
that in fact lie in the $76$ classes of maximal cyclic subgroups. There are in fact
$175$ classes of cyclic subgroups altogether, including the trivial group. We can
therefore take suitable powers of the $76$ words as a further set of $99$ words
defining elements of the group. 

First we label the $76$ words with the names given in \cite{Bclass}. These names will
later, of course, be identified with unions of conjugacy classes, but at this stage
they are simply names. We calculate the orders of the elements, and hence verify that the
numerical part of the name is indeed the order of the element. We define our
other $99$ words and their labels
as the obvious powers from the first line to the second line of each row of
Table~\ref{namedef}.
\begin{table}
\caption{\label{namedef}The names of our words}
$$
\begin{array}{cccccccccc}
\hline
70A & 66A & 60A & 60C & 56AB& 52A & 48B & 46AB & 44A & 66A\cr
35A & 33A & 30D & 30F & 28B & 26A & 24G & 23AB & 22B & 22A\cr\hline
42C & 40E & 40D & 60B & 60A & 40C & 38A & 36C & 36B & 34A\cr
21A & 20G & 20F & 20E & 20A & 20D & 19A & 18E & 18C & 17A\cr\hline
32CD & 32AB & 48B & 48A & 30F & 30A & 28E & 28A & 42A & 42B\cr
16D &  16C &  16B & 16A & 15B & 15A & 14E & 14D & 14A & 14B\cr\hline
42C & 26A & 24N & 24M & 24L & 24K & 24J & 24H & 24G & 24D\cr
14C & 13A & 12R & 12O & 12Q & 12M & 12J & 12F & 12G & 12D\cr\hline
36C & 36B & 36A & 60A & 60B & 22B & 20J & 20I & 20H & 20F\cr
12N & 12K & 12B & 12C & 12E & 11A & 10F & 10D & 10C & 10B\cr\hline
30A & 30E & 18F & 18E & 16H & 16G & 16F & 16E & 24J & 24M\cr
10A & 10E & 9B & 9A & 8M    & 8K  & 8H  & 8D  & 8J  & 8I\cr\hline
24I & 24K & 24C & 24B & 24A & 24E & 24N & 40D & 14D & 12T \cr
8G  & 8F  &  8E &  8C &  8B &  8A &  8N &  8L &  7A &  6K\cr\hline
12S &12R &12P &12O &12I &18A &30B &30A &30E &30C\cr
6J & 6I & 6H & 6G & 6C &  6D & 6A & 6B & 6E & 6F\cr\hline
12C & 10F &10B &8N & 8M & 8L & 8J & 8I & 8H & 8E\cr
4A  &  5B & 5A & 4J & 4H &4G & 4E & 4F & 4C & 4B\cr\hline
12E & 12T & 6K & 6A & 4J & 4I & 4A & 6A & 2B\cr
4D & 4I & 3B &  3A  & 2D & 2C & 2B & 2A & 1A\cr\hline
\end{array}
$$
\end{table}

At this stage, we have a list of $175$ words which give elements of the specified orders
in the Baby Monster. Our job now is to find invariants that
distinguish the alphabetical part of the name. 

\subsection{Invariants}
We compute only the invariants that 
\cite{Bclass} tells us are useful. Besides the order,
the  invariants we use for an element $x$ are of the following types:
\begin{itemize}
\item mod $2$ type: the trace $t_2(y)$ and the rank $r(y)$ of selected polynomials 
$y=p(x)$ in $x$, in the mod $2$ representation;
\item mod $3$ type: the trace $t_3(x^k)$ of selected powers of $x$, in the mod $3$ representation;
\item mod $5$ type: the trace $t_5(x)$ of $x$, in the mod $5$ representation.
\end{itemize}
The last two are expensive,
and are only used when we know they will in fact be useful.
\subsubsection{Odd-order elements}
We find $24$ cyclic subgroups of odd-order elements. For orders
$1$, $7$, $11$, $13$, $17$, $19$, $21$, $23$, $25$, $27$,
$31$, $33$, $35$, $39$, $47$ and $55$, the only
invariant we shall need is the order.
For the other orders, $3$, $5$, $9$ and $15$, the trace in the
$4370$ dimensional representation mod $2$ distinguishes two
names in each case:
$$\begin{array}{c|cccccccc|}
&3A&3B&5A&5B&9A&9B&15A&15B\cr\hline
t_2(x) &1&0&0&1&0&1&0&1
\cr\hline
\end{array}$$

\subsubsection{Elements of twice odd order}
For elements of order $38$, $46$, $66$ or $70$, no further
invariant is required.
In the remaining cases we compute the rank (or nullity) of $1+x$
on the $4370$ dimensional representation over the field of
order $2$. This  
turns out to be a sufficient 
invariant to distinguish all cases except the elements of
order $30$ and $42$.
The rank  
of $1+x$  
is tabulated below:
note that in the case $26A$ the rank is given incorrectly in \cite{Bclass} as $4196$
instead of $4198$.
$$
\begin{array}{ccccccccccccc}
\hline
2A & 2B & 2C & 2D & 6A & 6B & 6C & 6D & 6E & 6F \cr
1860 & 2048 & 2158 & 2168 & 3486 & 3510 & 3566 & 3534 & 3606 & 3604\cr
\hline
6G &6H & 6I & 6J & 6K & 10A & 10B & 10C & 10D & 10E\cr
3596 & 3610 & 3636 & 3638 & 3634 & 3860 & 3896 & 3918 & 3908 & 3920 \cr
\hline 10F & 14A &14B & 14C & 14D & 14E & 18A & 18B & 18C & 18D 
\cr
3932& 3996&4008 & 4048 & 4034 & 4052 & 4088 & 4090 & 4110 & 4124
\cr\hline
18E & 18F & 22A &22B & 26A & 26B & 30A/B & 30C & 30D & 30E 
\cr
4128 & 4122 & 4140 &4158 & 4198 & 4176 & 4190 & 4212 & 4206 & 4214 
\cr\hline
30F & 30GH & 34A & 34BC & 42A/B & 42C\cr
4224&4216&4238&4220&4242 & 4258\cr
\hline
\end{array}
$$
The cases $30A$ and $30B$ can be distinguished by the rank of $1+x^5$, which is $3510$ and $3486$
respectively. The cases $42A$ and $42B$ can be distinguished by the rank of $1+x^3$,
which is $3996$ and $4008$ respectively.

\subsubsection{Elements of order at least $28$} 
For elements of order $44$, $52$, and $56$, 
no further invariant is required. 
For elements of order $36$ and $60$, the rank of $1+x$ is
sufficient:
$$
\begin{array}{cccccc}
\hline
36A & 36B & 36C & 60A & 60B & 60C\cr
4226 & 4238 & 4248 & 4280 & 4286 & 4296
\cr\hline
\end{array}
$$
For the remaining element orders, $28$, $32$, $40$, and $48$, we have the following
values of the rank of $1+x$:
$$
\begin{array}{cccccccc}
\hline
28A/C & 28B/D & 28E & 32AB/CD & 40A/B/C & 40D & 40E &48A/B\cr
4188 & 4200 & 4210 & 4222 & 4242 & 4250 & 4258 & 4266
\cr\hline
\end{array}
$$
In particular, this invariant is of no help for elements of order $32$ or $48$.
All necessary cases can be separated by the trace mod $3$
of $x$ or $x^2$ or $x^7$:
$$
\begin{array}{c|cccccc}
&28A & 28C & 28B & 28D & 32AB & 32CD\cr\hline
t_3(x^2)&&&&&2&0\cr
t_3(x^7)&0&1&1&0\cr
\hline
\end{array}$$
$$\begin{array}{c|cccccc}
 & 40A & 40B & 40C & 48A & 48B\cr\hline
t_3(x)& 0&1&2 &0&1\cr
\hline
\end{array}
$$
\subsubsection{Elements of order $4$ and $8$}
The rank of $1+x$ distinguishes $7$ cases of elements of order $4$. Three of these split into two,
according to the trace on the $4371$-dimensional representation mod $3$.
An alternative invariant to distinguish $4H$ from $4J$ is the rank of
$(1+x)^3$ in the mod $2$ representation.
$$\begin{array}{c|cccccccccc|}
& 4A & 4B & 4C & 4D & 4E & 4F & 4G & 4H & 4I & 4J\cr\hline
r(1+x) & 3114 & 3114 & 3192 & 3192 & 3256 & 3202 & 3204 & 3266 & 3264 & 3266\cr
t_3(x) & 1 & 0 & 1 & 0 &&&& 0 && 1\cr
r(1+x)^3 &&&&&&&& 1082 && 1084
\cr\hline
\end{array}
$$
Similarly, for elements of order $8$, the rank of $1+x$ distinguishes $8$ cases,
one of which is split by the rank of $(1+x)^2$, while
the rank of $(1+x)^3$ splits two more:
$$
\begin{array}{ccccccccccc}
\hline
8A & 8B/C/E & 8D & 8F & 8H & 8G & 8I&8L & 8J & 8K/M & 8N\cr
3774 & 3738 & 3778 & 3780 & 3780 & 3810 & 3786 & 3786 & 3812 & 3818 & 3818\cr
&&&&&&3202&3204\cr
&&&2619 & 2620 &&&&& 2714 & 2717
\cr\hline
\end{array}
$$

The trace modulo $3$ distinguishes the remaining cases
$$
\begin{array}{c|ccccc|}
&8B & 8C & 8E & 8K & 8M\cr\hline
t_3(x) &1&0&2 & 1&2
\cr\hline
\end{array}$$

\subsubsection{Elements of order $12$ and $24$}
The rank of $1+x$ distinguishes $14$ cases of elements of order $12$:
$$
\begin{array}{c|ccccccc|}
& 12A/C/D & 12B & 12E & 12F & 12G/H & 12I & 12J\cr\hline
r(1+x) &3936 & 3942 & 3958 & 3996 & 3962 & 3964 & 3986\cr
\hline
\end{array}
$$
$$\begin{array}{c|ccccccc|}
&12K/M & 12L & 12N & 12O & 12P & 12Q/R/T & 12S\cr\hline
r(1+x) &3978 & 3966 & 4000 & 3982 & 3988 & 4002 & 4004\cr\hline
\end{array}
$$
All except $12A/D$ can be split using the trace mod $3$.
This last case would seem to require the trace mod $5$.
$$
\begin{array}{r|ccc|cc|cc|ccc|}
&12A & 12C & 12D & 12G & 12H & 12K & 12M & 12Q & 12R & 12T\cr\hline
t_3(x)&0&1&0&1&0 & 0&1 & 0&1&2\cr
t_5(x)&3&4&1&&&&&&&\cr\hline
\end{array}$$
Similarly for elements of order $24$, the rank of $1+x$
distinguishes $8$ cases
$$
\begin{array}{cccccccc}
\hline
24A/B/C/D & 24E/G & 24F & 24H & 24I/M & 24J & 24K & 24L/N\cr
4152 & 4164 & 4170 & 4182 & 4176 & 4178 & 4174 & 4186\cr
\hline
\end{array}
$$
Of these, we can distinguish $24I/M$ with the rank of $1+x^2$,
which is $3986$ and $3982$ respectively, and
$24E/G$ with the rank of $1+x^3$, which is $3774$ and $3778$
respectively.
All the rest are distinguished by the trace mod $3$,
apart from the case $24C/D$, which seems to require the trace mod $5$.
$$
\begin{array}{c|cccccc|}
&24A & 24B & 24C & 24D & 24L & 24N\cr\hline
t_3(x) & 1 & 0 & 2 & 2 & 1 & 2\cr
t_5(x) & 3 & 0 & 2 & 1 & & \cr
\hline
\end{array}
$$
\subsubsection{Elements of orders $16, 20$}
In these cases, the rank of $1+x$ distinguishes the following:
$$
\begin{array}{c|ccc|}
&16A/B & 16C/D/E/F & 16G/H \cr\hline
r(1+x) & 4072 & 4074 & 4094\cr
\hline
\end{array}$$
$$\begin{array}{c|ccccccc|}
&20A/B/C/D & 20E & 20F & 20G & 20H & 20I & 20J\cr\hline
r(1+x) & 4114 & 4128 & 4132 & 4148 & 4144 & 4138 & 4150\cr
\hline
\end{array}
$$
For the elements of order $20$, the rank of $1+x^2$ distinguishes $20B$
from $20A/C/D$, which are distinguished from each other by the
trace mod $3$:
$$\begin{array}{c|cccc|}
&20A&20B&20C&20D\cr\hline
r(1+x^2) & 3896 & 3908 & 3896 & 3896\cr
t_3(x) & 2 & 1 & 0 & 1\cr
\hline
\end{array}$$
For the elements of order $16$, the trace mod $3$
distinguishes $16A/B$, and separates $16C/E$ from $16D/F$.
Then $16C/E$ can be separated with the rank of $1+x^2$, and
$16G/H$ with the trace of $x^2$ mod $3$.
$$\begin{array}{c|ccccccc|}
& 16A & 16B & 16C & 16E & 16D/F & 16G & 16H\cr\hline
t_3(x) & 0 & 1 & 2 & 2 & 0 &  & \cr
t_3(x^2) & &&&&&1&2\cr
r(1+x^2) &&& 3780 & 3778&&&
\cr\hline
\end{array}
$$
There would appear to be no easily computed invariant which
distinguishes $16D$ from $16F$.

\subsection{Checking the approximate power maps}
We power up each of the given words, to every
relevant power (that is, every power dividing the element order), 
and compute the necessary invariants of the
resulting elements. 
We therefore know the power maps approximately.
In every case the power maps agree with the
character table in the {\Atlas} \cite{Atlas}. Indeed, some of the
power maps form part of the definition of our set of class
representatives, so the calculations in these cases
can in fact be omitted. This includes the class $16D$, which is
defined to be the square of the classes $32CD$. Hence it is
not necessary to find an invariant to distinguish $16D$ from $16F$,
in order to verify the power maps. All that remains in order to
verify that the power maps are actually correct, is, firstly, to
prove that the class list is correct, and secondly, to deal with
any issues concerning algebraically conjugate classes.
Details of the computations are given in \cite{ctblBM}.

\section{
Centralizers of prime order elements in the Baby Monster}%
\label{sect:oddcentralizers}

In this section
we determine the classes of prime order
elements, and the orders of their centralizers, in the Baby Monster. 
Much of this information comes from Stroth's 1976 paper \cite{Stroth}.
In cases where \cite{Stroth} does not give full information,
our strategy is first to use a certified copy of the Baby Monster from
\cite{webatlas} to give lower bounds on both the number of
conjugacy classes and the orders of the respective centralizers,
and then to use local arguments, together with information about the
permutation representation on the $\{3,4\}$-transpositions,
to show these are also upper bounds.
For technical reasons, we deal with the primes in
the order $2$, $3$, $7$, $17$, $11$, $13$, $19$, $23$, $5$, $47$, $31$.

\subsection{Fusion of involutions}
From \cite{Stroth} we see there are exactly four classes of involutions in $\B$,
with representatives labelled
$d$, $x_{36}(1)$, $d\beta$ and $dx_{33}(1)x_{32}(1)$ respectively.
In the {\Atlas} \cite{Atlas}, these are labelled $2A, 2B, 2C, 2D$ respectively.
The centralizers of $2A$ and $2B$ and $2C$ are given in detail in \cite{Stroth}.
The centralizer order of $2D$ is given together with a rough description of
the structure.

For the purposes of computation, it is necessary to match
these classes to the names given in Section~\ref{sect:classes}.
Note that any element of order $38$ powers into
class $2A$, and any element of order $22$ powers into class $2B$. An element of
order $34$ powers into either $2A$ or $2C$. We can now compute the rank of $1+x$
for involutions $x$ in the certified copy of the Baby Monster in dimension $4370$
mod $2$, obtaining the values $1860$ for class $2A$ and $2048$ for class $2B$.
For suitable $x$ obtained as the $17$th power of an element of order $34$, we obtain
the value $2158$, and for another involution we obtain $2168$, so these are
in class $2C$ and $2D$ respectively. 
Hence the names for involutions in Section~\ref{sect:classes}
are indeed the same as the {\Atlas} class names,
and we can now easily determine the class of
any explicitly given involution.

Following \cite{Stroth}, let $H = C_\B(d)\cong
2.{}^2E_6(2).2$ be a fixed $2A$ involution centralizer.
This group has ten classes of involutions, whose
labels in {\GAP} \cite{GAP} and the {\Atlas} are as follows:
$$\begin{array}{cccccccccc}
2&3&4&5&6&7&175&176&177&178\cr
2a&2b&2c&2d&2e&2f&2g&2h&2i&2j\cr
-1A & +2A & -2A & +2B & -2B & 2C & +2D & -2D & +2E & -2E
\end{array}$$
Since $H$ has shape $2.X.2$, it has an outer automorphism negating the outer classes,
and there is an arbitrary choice of which is which of classes $\pm 2D$.
For consistency with \cite{Stroth}, we choose $+2D$ to be the class
which fuses to $2A$ in the Baby Monster.
On the other hand, the classes $\pm 2A$ are distinguished in their
common centralizer $C=(2\times 2^{1+20}).\Uni_6(2).2$, in that
the involution in
the derived subgroup of $O_2(C)$ is in class $+2A$. It follows
that $+2A$ fuses to $2B$ in the Baby Monster, while $-2A$
fuses to $2A$ in the Baby Monster. 
Indeed, computation using explicit matrices, and suitable class invariants as above,
gives the full fusion of involutions from $H$ to $\B$. 
We find that {\GAP} classes $2,4,175$ fuse to $2A$, and classes $3,5$ fuse to $2B$,
while classes $176,177$ fuse to $2C$, and classes $6,7,178$ fuse to $2D$.

We will also need to know the fusion
of involutions from the subgroups $\Fi_{23}$ and $\HN$.
The easiest way to verify this is probably to use the words
in \cite{webatlas} to find these subgroups explicitly, and 
compute a suitable class invariant as described above.
We then see that
classes $2A, 2B, 2C$ in $\Fi_{23}$ fuse to
$\B$ classes $2A, 2B, 2D$ respectively, while
classes $2A, 2B$ in $\HN$ fuse to $2B,2D$ respectively in $\B$.

\subsection{The permutation representation on $\{3,4\}$-transpositions}
According to \cite{Stroth} the non-trivial suborbit lengths of $\B$ acting on the $13571955000$ cosets of $H$
are as follows: 
\begin{itemize}
\item $3^3.5.7.13.17.19=3968055$, with point stabilizer 
$2.2^{1+20}.\Uni_6(2).2$;
\item
$2^{12}.3^3.11.19=23113728$, with point stabilizer 
$2^2 \times F_4(2)$;
\item
$2^{20}.7.17.19=2370830336$, with point stabilizer
$\Fi_{22}.2$; and
\item
$2^8.3^3.5.7.11.13.17.19=11174042880$, with point stabilizer
 $2^{1+20}.\Uni_4(3).2^2$.
\end{itemize}

We now compute the permutation characters of the action of $H$ on 
the first three of these
suborbits. We use standard operations in {\GAP}, using only
the character tables of $H$ and certain of its subgroups.
For simplicity we use the {\GAP} labels for characters of $H$.

In the first case, the action on the suborbit 
is the permutation action of $^2E_6(2).2$ on the
cosets of the $\Uni_6(2)$ maximal parabolic, and is known to have rank $5$ (see, for example,
Theorem 4 in \cite{vasilyev}).
A straightforward combinatorial computation,
using {\GAP}, shows that the only way to get the character degrees adding to the correct number
is for the degrees to be $1+1938+48620+1828332+2089164$. The trivial character is
a constituent,
because it is a permutation character, 
leaving $16$ possibilities for the signs on the other four constituents.
It turns out that only one of these characters has non-negative values.
This character is the sum of the irreducibles labelled $1,3,5,13,15$ in {\GAP}.

In the second case, {\GAP} computes possible class fusions from $F_4(2)$ into $H$,
and we induce up the trivial character in each case. The answers are all the same.
The permutation character is a subcharacter of this induced character, 
and it is easy to determine
the character degrees, and then check all possibilities as above.
The answer is the sum of irreducibles numbered $1,5,17,24$.

In the third case, similarly, we compute possible class fusions from $\Fi_{22}{:}2$ into $H$.
There are then two possibilities for the induced trivial character, and they differ
by multiplying the outer elements of $\Fi_{22}{:}2$ by the central
involution of $H$. 
But we know that in the point stabilizer $\Fi_{22}{:}2$ the 
$2D$ involutions fuse to $2A$ in $\B$
(if necessary we can verify this computationally using the subgroup
$S_3\times \Fi_{22}{:}2$ in our certified copy of the Baby Monster), 
which distinguishes the two cases.
The answer is the sum of characters numbered 
$1,3,5,13,17,28,49,76,190,192,196,202,210,217$.

It is not necessary to compute the full permutation character
of $\B$ on the cosets of $H$, which would involve computing the 
fourth suborbit case as well.
Later on we will however need to compute the values on a few selected classes.

\subsection{Fusion of $3$-elements}
Computationally, using a certified copy
of the Baby Monster, and words provided in \cite{webatlas},
we find two subgroups $S_3\times \Fi_{22}{:}2$ and
$3^{1+8}{:}2^{1+6}.\Uni_4(2).2$, which normalize cyclic subgroups of order $3$.
The corresponding elements of order $3$ can be distinguished by the trace
in the $4370$ dimensional representation mod $2$, so do not fuse in $\B$.
We use the {\Atlas} labels $3A$ and $3B$ for these two conjugacy classes.

Conversely, note that
$\Fi_{23}$ contains a Sylow $3$-subgroup of $\B$ so every $3$-element in $\B$
is conjugate to an element of $\Fi_{23}$. Moreover, we know 
the fusion from $2.\Fi_{22}$ to $\Fi_{23}$,
and in particular, every $3$-class in $\Fi_{23}$ is represented in $2.\Fi_{22}$ and
therefore in $H$. 
Using the fact that $H$-classes 
$-2A$ and 
$+2D$ are in $2A$, 
and computing structure constants in $H$, we get that $H$-classes $3A$ and $3B$ fuse
in $\B$. 
Hence there are exactly two classes of 
elements of order $3$ in $\B$. 

We now show that a $3A$-element $x$ has centralizer $C_{\B}(x)\cong 3\times \Fi_{22}{:}2$ in $\B$.
We know its centralizer is at least that (either computationally,
as above, or see \cite{Stroth}). On the other hand,
the number of $\B$-conjugates of $x$ is at least one-third of
the product of the length of the whole orbit with the length of the relevant
suborbit. 
This number $13571955000\times 2370830336/3$ is equal to the index of $3\times \Fi_{22}{:}2$ in $\B$, and the claim
is proved.

Next we show that the subgroup $3^{1+8}.2^{1+6}.\Uni_4(2)$ computed above is the full centralizer
of a $3B$-element. To do this we need to know the value on $3B$ of the full
permutation character of $\B/H$. Equivalently, the value 
of the permutation character of the last orbit above on $H$-class $3C$.
Recall that the point stabilizer in the last orbit is $H_5=2^{1+20}.\Uni_4(3).2^2$.
We use the {\GAP} function {\tt PossibleClassFusions} applied to the character tables of
$\Uni_4(3)$ and $2.\Uni_4(3)$ to get
 the 
 fusion of $3$-elements from $H_5$ 
to $H$.
The result is that {\Atlas} class $3A$ in $\Uni_4(3)$ fuses to $3C$ in $H$, 
while all other classes of elements of order $3$ fuse to $3A$ or $3B$ in $H$. 
Hence the value on $H$-class $3C$ of the permutation character of $H$ on this orbit 
is $|C_H(3C)|/|C_{H_5}(3A)|=2^{11}.3^9/2^8.3^6=2^3.3^3=216$.
Therefore the value of the whole permutation character of $B/H$ on class $3C$ is $1620$.

Hence we know $|C_\B(3B)|/|C_H(3C)|=1620=2^2.3^4.5$, 
and the $3C$-centralizer in $H$ has order $2^{11}.3^9$,
so we deduce the order of the
$3B$-centralizer in $\B$ is $2^{13}.3^{13}.5$ and the claim follows.

\subsection{
Elements of orders $7,11,13,17,19$ and $23$}
From \cite{Stroth} (Lemma 6.11)
we get $|N_\B(7^2)|=2^6.3^2.7^2$ and
$|C_\B(7^2)|=2^2.7^2$. The only subgroup of $\GL_2(7)$ of order $2^4.3^2$ is
$3\times 2S_4$, which is transitive on non-zero elements of $7^2$. Hence there
is a single class of elements of order $7$. (This can also be verified
computationally in a certified copy of the Baby Monster.)

In $H$ we have a $7A$-centralizer $7 \times 2.\Lin_3(4).2$.
To show that the centralizer in $\B$ is no bigger, we follow the same strategy
as for $3B$ elements above, although it is slightly more
complicated since both classes $7A$ and $7B$ fuse to
$7A$ in $\B$. The $7$-elements in $H_5$ fuse to class $7B$ in $H$.
Hence the value of the permutation character of $\B/H$ on $H$-class $7A$
is $121$.
(As a check,
$|C_{H}(7B)|=4704=2^5.3.7^2$ and $|C_{H_5}(7B)|=2^4.7$, so the character value
of the last orbit on $H$-class $7B$ is $2.3.7=42$. This implies the value of the
permutation character of $\B/H$
on $H$-class $7B$ is also $121$, as it must be.)
Therefore
$|C_\B(7A)|/|C_H(7A)| + |C_\B(7A)|/|C_H(7B)|=121$, so that $|C_\B(7A)|=2^8.3^2.5.7$.

For the remaining primes in the list, $11$, $13$, $17$, $19$ and $23$,
most of the information we need is already
in \cite{Stroth}.
Lemma 6.13 of \cite{Stroth} says the centralizer order of an
element of order $17$ is $2^2.17$,
and the normalizer has order $2^6.17$, so there is a single class of elements of
order $17$.
Lemma 6.8 of \cite{Stroth} says that
the order of the Sylow $11$-normalizer is $2^4.3.5^2.11$, and the centralizer of
an element of order $11$ is
$S_5\times 11$. Hence the normalizer is $S_5\times 11{:}10$.
In Lemma 6.12 of \cite{Stroth} there are two possibilities for the 
normalizer of an element of order $13$. 
But the normalizer of such an element in $F_4(2)$ is just $13{:}12$,
so from the proof of Lemma 6.12 we get that  the $13$-centralizer in $\B$
is $13\times S_4$, and the normalizer is $13{:}12\times S_4$.

The Sylow $19$-subgroup is self-centralizing in $H/\langle d\rangle$,  so the Sylow $2$-subgroup of $C_\B(19)$
is of order $2$, containing a $2A$-element. Since $|\Fi_{22}|$ 
is not divisible by $19$, that forces $N_\B(19)$ to lie in $H$.
Lemma 6.20 of \cite{Stroth} says that
the $23$-normalizer contains $2\times 23{:}11$, and 
that the Sylow $2$-subgroup of the $23$-normalizer has order $2$;
we know (from Lemmas 7.13, 7.14, 7.15 and 7.17 of \cite{Stroth}) that all Sylow subgroups
of the normalizer are cyclic. From the discussion earlier in this section,
we know that the normalizer does not contain elements
of order $7, 13,17$, or $19$. The normalizing $11$ rules out $47$ and $31$, by the Frattini argument.
This leaves $3,5$. We know the $3$-centralizers, so $3$ is ruled out.
Finally $5$ is ruled out because $|\B|/(2.5.11.23)\not\equiv1\pmod {23}$.

\subsection{The elements of order $5$}

The subgroup $\HN$ (constructed explicitly in our certified copy of $\B$) 
contains a full Sylow $5$-subgroup.
Every element of order $5$ in $\HN$ centralizes an involution, which we know
fuses to $2B$ or $2D$ in $\B$. Moreover, every element of order $5$ in $C_\B(2B)$
or $C_\B(2D)$ centralizes an element of order $3$.
But $C_\B(3A)$ and $C_\B(3B)$ contain just one class of elements of order $5$ each,
so there are at most two classes of elements of order $5$ in $\B$.
On the other hand, we find two classes of $5$-elements with different traces.
Hence there are exactly two classes.

The usual argument gives the order of $C_\B(5A)$. We have
$|C_H(5A)|=2^8.3^2.5^2.7$ and $|C_{H_5}(5A)|=2^7.5$ so the value of the 
permutation character of the last orbit on this class
is $2.3^2.5.7=630$.
Hence the full permutation character has value $630+470=1100$ on
$\B$ class $5A$. Therefore $|C_\B(5A)|=1100|C_H(5A)|=2^{10}.3^2.5^4.7.11$,
which is the order of $5\times \HS{:}2$.
Hence the $5A$-normalizer is $5{:}4 \times \HS{:}2$.

Computationally, using the matrices and words provided in \cite{webatlas}, 
we find a subgroup $5^{1+4}.2^{1+4}.A_5.4$, normalizing
a cyclic group of order $5$, which must therefore be of $5B$ type.
We shall show that the normalizer is no bigger than this.
We know that there is no $5B$ element in the centralizer of any element
of order $7,11,13,17,19,23$, or of a $3A$. Also $47$ and $31$ are $3\mod 4$
and do not centralize an involution, so do not centralize a $5B$ by the Frattini argument.
Hence the centralizer of a $5B$ is a $\{2,3,5\}$-group, and contains the
full Sylow $5$-subgroup of $\B$, so 
only the Sylow $2$- or $3$-subgroup could grow.

Now the centralizer of a $5$-element in $C_\B(3B)=3^{1+8}.2^{1+6}.\Uni_4(2)$ is
just a cyclic group of order $30$, so the Sylow $3$-subgroup of the
$5B$-centralizer has order $3$.
Since $C_\B(2A)$ and $C_\B(2C)$ contain no $5B$, we look in $C_\B(2B)$ and $C_\B(2D)$.
In $C_\B(2B)$ only the $5A$ class of $\Co_2$ fuses to $\B$ class $5B$, and we
see the centralizer $(2^{1+2}\times 5^{1+2}).2A_4$ of order $2^6.3.5^3$.
In $C_\B(2D)$ we see centralizer order $2^7.3.5^2$. In neither case does the
Sylow $2$-subgroup grow.
Thus we know the orders of all the Sylow subgroups of $C_\B(5B)$,
and therefore the order of $C_\B(5B)$.

\subsection{Primes $47$ and $31$}
The order of the $47$-normalizer now divides $47.2.23.31$ and Sylow's theorem
implies it is $47.23$. Finally the order of the $31$-normalizer divides
$31.30$, so is $31.15$ by Sylow's Theorem.

\section{Obtaining the class list}
\label{classlist}
Our strategy for obtaining the list of conjugacy classes in the Baby Monster
is first to determine the classes of even order elements, by computing the
character tables of subgroups containing the 
four distinct involution centralizers, and noting down
the conjugacy classes of elements in each subgroup that power to the
relevant involution class. (The centralizers of involutions in classes
$2A$, $2B$, $2C$ are in fact maximal, although it is not
necessary to know this, so we have no choice but to use the involution
centralizer itself in these cases.) At the same time, we note down the length of each such
class. A similar computation for odd-order elements in the centralizers
of elements of odd prime order is trivial in comparison.

In fact, there is a great deal of redundancy in the information that we have computed,
and classes of elements whose order is divisible by two primes can be computed
in two different ways. This provides a robust check on these results, in particular
for the large number of classes of elements of order divisible both by $2$ and by an
odd prime.

\subsection{Involution centralizers in the Baby Monster}%
\label{sect:invcentralizers}

The character table of the $2A$-centralizer is known
by \cite{BMO17} and the computations shown in~\cite{ctblatlas}.
The 
$2C$-centralizer has the structure
$(2^2 \times F_4(2)).2 < D_8 \times F_4(2).2$, and
its character table is determined by those of the subgroups
$2^2$ and $F_4(2)$ and the factor groups $D_8$ and $F_4(2).2$,
hence it is known.

The 
$2D$-centralizer is 
contained in subgroups of the
structure $2^{(8+1)+16}.\Sp_8(2)$ in $\B$.
Such subgroups can be constructed explicitly in a certified copy
of $\B$, 
using the straight line program from~\cite{webatlas}.
The character table of this subgroup can be verified
by restricting the $2$-modular degree $4370$ representation of $\B$
to the subgroup,
finding a faithful $180$-dimensional
subquotient of this module,
and computing the character table from this matrix representation
using the \textsf{MAGMA} computer algebra system~\cite{MAGMA}.
(This had been done by E.~O'Brien in 2007,
but we repeat the computations in order to make sure
that only 
explicitly verified data are used.)
In particular, this verification includes a verification
that the given subgroup contains the full $2D$-centralizer.

The 
$2B$-centralizer has the structure $2^{1+22}.\Co_2$.
The character table has been computed in~\cite{Pah07}
but the arguments assume the character table of $\B$.
In the remainder of this section,
we describe briefly how we verify this character table.
Full details can be found in \cite{ctblBM}.

First we restrict the certified
$3$- and $5$-modular representations of degree $4371$
of $\B$ to the
$2B$-centralizer,
using the straight line program from~\cite{webatlas};
the composition factors of the module have the dimensions $23$, $2300$,
and $2048$ in both cases.
Next, we find an orbit of length $4600$ in the $2300$-dimensional module
over the field of order $3$.
The action on this orbit yields a faithful permutation representation of the
factor group $2^{22}.\Co_2$.
We compute class representatives for this factor group,
and let {\MAGMA} compute its character table.

The $2048$-dimensional module is faithful.
We compute the class fusion under the epimorphism from $2^{1+22}.\Co_2$
to $2^{22}.\Co_2$, and the Brauer characters of our $3$- and $5$-modular
representations for this module.
Now $2^{1+22}$ has a unique faithful irreducible representation in every
characteristic except $2$, and this representation has dimension $2^{11}=2048$,
and extends uniquely to $2^{1+22}.\Co_2$. If $\chi$ denotes the character
of the ordinary representation of $2^{1+22}.\Co_2$ obtained in this way,
then
all faithful irreducible characters of $2^{1+22}.\Co_2$ arise as tensor
products of $\chi$ with the irreducibles of the factor group $\Co_2$.
In particular the Brauer characters computed above lift to $\chi$,
and therefore we obtain the values of $\chi$ on all classes of elements
whose order is not divisible by $15$. For details of how the
remaining values were obtained, see \cite{ctblBM}.
 
Once the character tables of (overgroups of) all the involution
centralizers are available, we can read off from these tables all the conjugacy
classes of elements that power to each of the involutions, together with the
centralizer orders. This gives us a complete list of all conjugacy classes
of even-order elements.

\subsection{Elements of odd order} 
For primes $p\ge 11$ it is now almost a triviality to write down the classes
of elements of odd order divisible by $p$. 
For $p=7$, we have that $N_\B(7A)$ is contained in $H$, so the 
relevant classes can be read
off from the classes of elements of $H$ that power into class $7A$. 
(Note that the $7A$-centralizer in $H$ has the
shape $7\times 2.\Lin_3(4){:}2_2$, and the $2_2$ automorphism swaps the $\Lin_3(4)$-classes
$5A$ with $5B$.)

For elements powering into $5A$, we read off the classes 
and their centralizer orders from the
{\Atlas} character table of $\HS{:}2$.
Similarly, for elements powering into $3A$, use the table for $\Fi_{22}{:}2$,
but note that there are some classes missing in the {\Atlas} character table for
$\Fi_{22}{:}2$: these only affect the calculations for elements of order $30$,
which have already been dealt with in the $5A$-centralizer.

In the cases $3B$ and $5B$ again, {\GAP} contains character tables of the
respective normalizers. However, it is not recorded exactly what information
was used to calculate these tables. Therefore we re-calculate them (see \cite{ctblBM}).
In conclusion,
we find that the list of odd-order elements and their centralizer orders
agrees with the {\Atlas}.

\section{Computing the irreducible characters of the Baby Monster}
\label{sect:irreducibles}

From the previous sections, we know that $\B$ contains subgroups
of the structures $2.{}^2E_6(2).2$, $\Fi_{23}$, and $\HN.2$.
The ordinary character tables of these groups have been verified
(see~\cite{BMO17}) and thus may be used in our computations.
The class fusions from these subgroups to $\B$ can be computed
with the methods available in {\GAP} \cite{GAP}.
Moreover, in Section~\ref{sect:invcentralizers},
we have computed the character table of the 
$2B$-centralizer in $\B$. 
The class fusion from $2^{1+22}.\Co_2$ to $\B$ is determined
by evaluating the three representations of $\B$ at the class representatives
of $2^{1+22}.\Co_2$,
and applying the invariants from Section~\ref{sect:classes}.

Thus we can induce the irreducible characters from these subgroups
to $\B$.
Using the power maps of $\B$,
we induce also the irreducibles of all cyclic subgroups of $\B$.
Now we proceed in two steps.

In the first step,
we assume that $\B$ has an ordinary irreducible representation
$\chi$ of degree $4371$ such that the reductions modulo $3$ and $5$
are (irreducible and) equivalent to the representations we have used
in the previous sections,
and such that the reduction modulo $2$
has one trivial composition factor and one that is equivalent to
the representation we have used above.
Then the Brauer character values of our representations yield the values
of $\chi$, except on the classes of elements with order divisible by $30$,
and the missing values are uniquely determined by the obvious bounds.
If we add $\chi$ and the trivial character of $\B$ to the list of
induced characters then applying standard character-theoretic techniques
such as LLL reduction yields
a complete list of irreducible characters for $\B$,
which coincides with the characters in the {\Atlas} table of $\B$.

In the second step,
we do not want to assume the existence of the ordinary character $\chi$,
and try to apply the character-theoretic criteria to the safe list of
induced characters.
This way, we do not get any irreducible character.
However, we can show that $30$ 
class functions from the list of irreducibles
computed in the first step lie in the lattice spanned by the induced
characters.
Thus these 
class functions are verified as irreducible characters of $\B$.
Now we form symmetrizations and tensor products of the known
irreducible characters,
and the lattice spanned by the known characters of $\B$ contains
all the missing irreducibles computed in the first step.
Thus we are done.
Again, the details of these constructions can be found in \cite{ctblBM}.

\section*{Acknowledgements}
We thank Chris Parker for significant contributions to the original version of this paper,
and we thank the referee for helpful comments that enabled us to avoid the need for them.


\end{document}